\documentclass[12pt,reqno]{article}
\usepackage{authblk}
\usepackage{hyperref}
\usepackage{amssymb}
\usepackage{amsthm}
\usepackage{amsmath}
\usepackage{multirow}
\usepackage{float}
\usepackage{graphicx}
\usepackage{listings}
\usepackage{MnSymbol}
\usepackage{nicematrix}
\usepackage{pdfpages}
\newtheorem{theorem}{Theorem}[section]
\newtheorem{corollary}[theorem]{Corollary}
\newtheorem{lemma}[theorem]{Lemma}

\newtheorem{definition}[theorem]{Definition}

\newtheorem{remark}[theorem]{Remark}

\DeclareMathOperator{\vect}{vec}
\DeclareMathOperator{\diag}{diag}
\begin{document}
\title{\Large\bfseries Efficient tensor-based approach to solving linear systems involving Kronecker sum of matrices}
\author[1]{Ahmad Y. Al-Dweik}
\author[2]{Abdallah Sayyed-Ahmad}
\affil[1]{Department of Mathematics, Birzeit University, Ramallah, Palestine; aaldweik@birzeit.edu}
\affil[2]{Department of Physics, Birzeit University, Ramallah, Palestine;
asayyeda@birzeit.edu}
\maketitle
\begin{abstract}
A novel tensor-based formula for solving the linear systems involving Kronecker sum is proposed.  Such systems are directly related to the matrix and tensor forms of Sylvester equation. The new tensor-based formula demonstrates the well-known fact that a Sylvester tensor equation has a unique solution if the sum of spectra of the matrices does not contain zero. We have showcased the effectiveness of the method by efficiently solving the 2D and 3D discretized Poisson equations, as well as the 2D steady-state convection-diffusion equation, on a rectangular domain with Dirichlet boundary conditions. The results suggest that this approach is well-suited for high-dimensional problems.
\end{abstract}
\section{Introduction}
A \textit{Sylvester equation}\cite{Sylvester1884} refers to a matrix equation of the form $$\mathbf{A}\mathbf{X}+\mathbf{X}\mathbf{B}=\mathbf{Y},$$ which is traditionally solved using the Bartels–Stewart algorithm\cite{BartelsStewart1972} in which $\mathbf{A}$ and $\mathbf{B}$ are transformed to the \textit{Schur} form. Subsequently, the resulting triangular system is solved via back-substitution. This algorithm has a computational complexity $\mathcal {O}(n^{3})$. A Sylvester equation has a unique solution for $\mathbf{X}$ when the sum of spectra of the matrices $\mathbf{A}$ and $\mathbf{B}$ does not contain zero (i.e. $\lambda_1+\lambda_2 \ne 0,~\forall \lambda_1 \in \sigma(\mathbf{A}),\lambda_2 \in \sigma(\mathbf{B})$).

The third order \textit{Sylvester tensor equation} $$\mathcal{X} \times_1 \mathbf{A}+\mathcal{X} \times_2 \mathbf{B}+\mathcal{X} \times_3 \mathbf{C}=\mathcal{Y}$$ is widely used in numerous scientific and engineering applications, including finite difference \cite{Bai2003}, thermal radiation \cite{Li2010}, information retrieval \cite{Li2015}, and finite elements \cite{Grasedyck2004}.
In the past decade, many theoretical and numerical methods for the solution of the Sylvester tensor equation have been developed\cite{Chen2012,Chen2013,Ballani2013,Panjeh2015,Ali2016, Zhang2021,Zhang2022,Chen2023}. A Sylvester tensor equation has a unique solution if and only if the sum of spectra of the matrices $\mathbf{A}, \mathbf{B}$ and $ \mathbf{C}$ does not contain zero \cite{Chen2013}. The linear system that involves the Kronecker sum of matrices is connected to the matrix and tensor forms of the Sylvester equation, as will be detailed in Section 3.

The $2D$ discrete Laplacian on a regular grid with homogeneous Dirichlet boundary conditions can be expressed as $\mathbf{D}_{xx} \oplus \mathbf{D}_{yy}$, where $\mathbf{D}_{xx}$ and $\mathbf{D}_{yy}$ represent the $1D$ discrete Laplacians in the $x$- and $y$-directions, respectively. Both \(\mathbf{D}_{xx}\) and \(\mathbf{D}_{yy}\) must represent the discrete Laplacian with homogeneous Dirichlet boundary conditions at the end points of the \(x\)- and \(y\)-intervals. This ensures that the resulting 2D discrete Laplacian correctly reflects the homogeneous Dirichlet boundary conditions across the entire boundary of the rectangular domain \cite{Lyche2020}. Similarly, the 3D discrete Laplacian on a regular grid with homogeneous Dirichlet boundary conditions can be expressed as \(\mathbf{D}_{xx} \oplus \mathbf{D}_{yy} \oplus \mathbf{D}_{zz}\), where \(\mathbf{D}_{xx}\), \(\mathbf{D}_{yy}\), and \(\mathbf{D}_{zz}\) are the $1D$ discrete Laplacians in each of the three directions. Interestingly, this formulation can be generalized to \(N\) dimensions \cite{Ballani2013}.  

To our knowledge, no direct solver currently exists for linear systems involving the Kronecker sum of matrices that accounts for the Kronecker sum structure. This gap motivates us to develop a general formula for solving such linear systems.

The outline of the paper is as follows. In the next section, we review key definitions and foundational facts. Subsequently, we present the relation between the linear system involving the Kronecker sum of matrices or tensors in the form of the Sylvester equation in section 3. In section 4, we propose a new tensor-based formula for the solution of the linear system involving the Kronecker sum of matrices. Finally, we investigate the performance and accuracy of the solutions of the discrete Poisson and steady state convection diffusion equations in section 5.
\section{Preliminaries}
In this section, we will present essential definitions and results that will be utilized in the subsequent discussion. For more detailed information, we refer the reader to the references \cite{Kolda2006, Kolda2009, Qi2017}.
\begin{definition}(Vec-operation). 
For any $\mathbf{A} \in \mathbb{C}^{m \times n }$, we define the vector
$$\vect(\mathbf{A})=[a_{11},\dots,a_{m1},a_{12},\dots,a_{m2},\dots., a_{1n},\dots,a_{mn}]^T,$$
where the columns of $\mathbf{A}$ are stacked on top of each other.
\end{definition}
\begin{definition}(Kronecker product). 
Let $\mathbf{A} \in \mathbb{C}^{r \times s}$ and $\mathbf{B} \in \mathbb{C}^{m\times n}$ be two
matrices. Then the $rm \times sn$ matrix
$$
{\mathbf{A}} \otimes {\mathbf{B}} = \left[ {\begin{array}{*{20}c}
   {a_{11} {\mathbf{B}}} &  \cdots  & {a_{1s} {\mathbf{B}}}  \\
    \vdots  &  \ddots  &  \vdots   \\
   {a_{r1} {\mathbf{B}}} &  \cdots  & {a_{rs} {\mathbf{B}}}  \\
\end{array}} \right],
$$
is the Kronecker product of $\mathbf{A}$ and $\mathbf{B}$.
\end{definition}
\begin{definition}(Kronecker sum). 
If  $\mathbf{A} \in \mathbb{C}^{n \times n}$, $\mathbf{B} \in \mathbb{C}^{m\times m}$ and $\mathbf{I}_k$ denotes the $k \times k$ identity matrix, then we can define the Kronecker sum, $\oplus$, by
$$
{\mathbf{A}} \oplus {\mathbf{B}} = {\mathbf{A}} \otimes {\mathbf{I}}_m  + {\mathbf{I}}_n  \otimes {\mathbf{B}}.
$$
\end{definition}
\begin{definition}(Hadamard product). 
If  $\mathbf{A}$ and $\mathbf{B} \in \mathbb{C}^{n \times m}$, the Hadamard product 
$\mathbf{A} \odot \mathbf{B}$ is a matrix of the same dimension as the operands, with elements given by
$$
(\mathbf{A} \odot \mathbf{B})_{ij}  = a_{ij} b_{ij} .
$$
\end{definition}
\begin{definition}(Outer product of tensors). 
If ${\cal T} \in \mathbb{C}^{k_1  \times k_2 \times \dots \times k_m}$ and ${\cal X} \in \mathbb{C}^{l_1  \times l_2 \times \dots \times l_n}$ are tensors, then their outer product ${\cal T} \circ {\cal X}$ is a tensor  $\in \mathbb{C}^{k_1  \times k_2 \times \dots \times k_m \times l_1  \times l_2 \times \dots \times l_n}$ with entries
$$
({\cal T} \circ {\cal X})_{i_1,i_2,\dots,i_m,j_1,j_2,\dots,j_n}={\cal T}_{i_1,i_2,\dots,i_m}{\cal X}_{j_1,j_2,\dots,j_n}.
$$
\end{definition}
\begin{definition}(The $n$-mode product of a tensor with a matrix). 
The $n$-mode product of a tensor ${\cal T} \in \mathbb{C}^{d_1  \times d_2 \times \dots \times d_N}$ with a matrix  $\mathbf{V} \in \mathbb{C}^{k \times d_n}$ is denoted by
${\cal T}  \times_n \mathbf{V}$ and is of size $d_1  \times \dots \times d_{n-1}\times k \times d_{n+1}\dots \times d_N.$
Elementwise, we have
$$
({\cal T} \times _n \mathbf{V})_{ j_1,\dots,j_{n-1},i, j_{n+1}, \dots,j_N} = \sum\limits_{j_n  = 1}^{d_n } {\cal T}_{j_1 ,j_2 ,\dots, j_N }\mathbf{V}^{}_ {i,j_n}.
$$
\end{definition}
\begin{remark}(Compatibility of Vectorization with Kronecker products).
$$\vect(\mathbf{A}\mathbf{B}\mathbf{C}) = (\mathbf{C}^T  \otimes \mathbf{A})\vect(\mathbf{B})$$
for matrices $\mathbf{A}$, $\mathbf{B}$, and $\mathbf{C}$ of dimensions $k\times l$, $l\times m$, and $m\times n$ respectively. 
\end{remark}
\begin{remark}(Compatibility of Vectorization with Hadamard products).
$$
\vect (\mathbf{A} \odot \mathbf{B})=\vect (\mathbf{A}) \odot \vect(\mathbf{B})
$$
for matrices $\mathbf{A}$ and $\mathbf{B}$ of dimensions $n\times m$.
\end{remark}
\begin{remark}(Connection of Outer product with the Kronecker product).
In the case of column vectors, the Kronecker product can be viewed as a form of vectorization of the outer product.
$$
\vect (\mathbf{u} \circ \mathbf{v}) =\mathbf{v} \otimes \mathbf{u},
$$
for vectors $\mathbf{u}$ and $\mathbf{v}$. Note that the vectorization transforms the outer product of vectors into the kronecker proudct but reverses the order of the vectors.
\end{remark}
\section{Transform the linear system involving Kronecker sum of matrices into matrix and tensor form}
 Here, we explore the relationship between the linear system considered in two dimensions and the well-known Sylvester equation. Additionally, we extend this connection to higher dimensions by relating the linear system to a Sylvester tensor equation.
\begin{lemma}\label {lemma1}(Vectorization with $n$-mode product of a 3rd-order tensor with matrices).
Given a tensor ${\cal T} \in \mathbb{C}^{d_1  \times d_2  \times d_3 }$ and matrices $\mathbf{A},\mathbf{B}$ and $\mathbf{C}$ in $\mathbb{C}^{k_1  \times d_1 } ,\mathbb{C}^{k_2  \times d_2 } ,\mathbb{C}^{k_3  \times d_3 },$  respectively. Then
$$\vect\left({\cal T} \times_1 \mathbf{A} \times_2 \mathbf{B} \times_3 \mathbf{C}\right)=\left(\mathbf{C}\otimes \mathbf{B} \otimes \mathbf{A}\right) \vect({\cal T})$$
\end{lemma}
\proof
We can express a tensor ${\cal T} \in \mathbb{C}^{d_1  \times d_2  \times d_3 }$ in term of the basis vectors
\begin{equation}
{\cal T}= \sum\limits_{i_1  = 1}^{d_1 }\sum\limits_{i_2  = 1}^{d_2 }\sum\limits_{i_3  = 1}^{d_3 } {\cal T}_{i_1 ,i_2 ,i_3 }(\mathbf{e}_{i_1}^{d_1} \circ \mathbf{e}_{i_2}^{d_2} \circ \mathbf{e}_{i_3}^{d_3}),\\
\end{equation}
where $\mathbf{e}_{i}^{n}$ is the $i$-th canonical vector of $\mathbb{C}^{n}$. Then the modal product of a 3rd-order tensor ${\cal T}$ with matrices $\mathbf{A},\mathbf{B}$ and $\mathbf{C}$ can be written as
\begin{equation} 
{\cal T} \times_1 \mathbf{A} \times_2 \mathbf{B} \times_3 \mathbf{C}=\sum\limits_{i_1  = 1}^{d_1 }\sum\limits_{i_2  = 1}^{d_2 }\sum\limits_{i_3  = 1}^{d_3 } {\cal T}_{i_1 ,i_2 ,i_3 }(\mathbf{A} \mathbf{e}_{i_1}^{d_1} \circ \mathbf{B} \mathbf{e}_{i_2}^{d_2} \circ \mathbf{C} \mathbf{e}_{i_3}^{d_3}),\\
\end{equation}

\begin{equation} 
\begin{array}{rl}
\vect\left({\cal T} \times_1 \mathbf{A} \times_2 \mathbf{B} \times_3 \mathbf{C}\right)&=\sum\limits_{i_1  = 1}^{d_1 }\sum\limits_{i_2  = 1}^{d_2 }\sum\limits_{i_3  = 1}^{d_3 } {\cal T}_{i_1 ,i_2 ,i_3 }\vect( \mathbf{A} e_{i_1}^{d_1} \circ \mathbf{B} e_{i_2}^{d_2} \circ  \mathbf{C} e_{i_3}^{d_3}),\\
&=\sum\limits_{i_1  = 1}^{d_1 }\sum\limits_{i_2  = 1}^{d_2 }\sum\limits_{i_3  = 1}^{d_3 } {\cal T}_{i_1 ,i_2 ,i_3 }( \mathbf{C} \mathbf{e}_{i_3}^{d_3} \otimes \mathbf{B} \mathbf{e}_{i_2}^{d_2} \otimes  \mathbf{A} \mathbf{e}_{i_1}^{d_1}),\\
&=( \mathbf{C} \otimes \mathbf{B} \otimes  \mathbf{A} )\sum\limits_{i_1  = 1}^{d_1 }\sum\limits_{i_2  = 1}^{d_2 }\sum\limits_{i_3  = 1}^{d_3 } {\cal T}_{i_1 ,i_2 ,i_3 }(\mathbf{e}_{i_3}^{d_3} \otimes  \mathbf{e}_{i_2}^{d_2} \otimes \mathbf{e}_{i_1}^{d_1}),\\
&=( \mathbf{C} \otimes \mathbf{B} \otimes  \mathbf{A} )\sum\limits_{i_1  = 1}^{d_1 }\sum\limits_{i_2  = 1}^{d_2 }\sum\limits_{i_3  = 1}^{d_3 } {\cal T}_{i_1 ,i_2 ,i_3 }\vect(\mathbf{e}_{i_1}^{d_1} \circ  \mathbf{e}_{i_2}^{d_2} \circ \mathbf{e}_{i_3}^{d_3}),\\
&=( \mathbf{C} \otimes \mathbf{B} \otimes  \mathbf{A} )\vect({\cal T} ).\\
\end{array}
\end{equation}
\endproof
\subsection{Sylvester equation}
Let  $ \mathbf{A},  \mathbf{B},  \mathbf{X}$ and $\mathbf{Y}$ $\in \mathbb{C}^{n \times n}$ be square matrices such that $\mathbf{x}=\vect( \mathbf{X})$ and $\mathbf{y}=\vect(\mathbf{Y})$. Then
\begin{equation}
\begin{array}{rl}
\left( \mathbf{A} \oplus  \mathbf{B} \right)\mathbf{x}=\mathbf{y}&\iff \left( \mathbf{A} \otimes  \mathbf{I} + \mathbf{I} \otimes  \mathbf{B}\right)\vect( \mathbf{X})=\vect( \mathbf{Y})\\
&\iff \left(\mathbf{A} \otimes\mathbf{I} \right)\vect(\mathbf{X})+\left(\mathbf{I} \otimes \mathbf{B}\right)\vect(\mathbf{X})=\vect(\mathbf{Y})\\
&\iff \vect\left(\mathbf{I}\mathbf{X}\mathbf{A}^T\right)+\vect\left(\mathbf{B}\mathbf{X}\mathbf{I}\right)=\vect(\mathbf{Y})\\
&\iff 
\vect\left(\mathbf{X}\mathbf{A}^T\right)+\vect\left(\mathbf{B}\mathbf{X}\right)=\vect(\mathbf{Y})\\
&\iff 
\mathbf{X}\mathbf{A}^T+\mathbf{B}\mathbf{X}=\mathbf{Y}.\\
\end{array}
\end{equation}
\subsection{Sylvester tensor equation}
Let $ \mathbf{A},  \mathbf{B}$ and $\mathbf{C}$ in $\mathbb{C}^{n_1  \times n_1 } ,\mathbb{C}^{n_2  \times n_2 } ,\mathbb{C}^{n_3  \times n_3 },$  respectively, and $\mathcal{X}$ and $\mathcal{Y}$ $\in \mathbb{C}^{n_1 \times n_2 \times n_3}$ be tensors such that $\mathbf{x}=\vect(\mathcal{X})$ and $\mathbf{y}=\vect(\mathcal{Y})$. Then
\begin{equation}
\begin{array}{rl}
\left(\mathbf{A} \oplus \mathbf{B} \oplus \mathbf{C} \right)\mathbf{x}=\mathbf{y}&\iff \left(\mathbf{A} \otimes \mathbf{I} \otimes \mathbf{I}  +\mathbf{I} \otimes \mathbf{B}\otimes \mathbf{I}+\mathbf{I} \otimes \mathbf{I}\otimes \mathbf{C} \right)\vect(\mathcal{X})=\vect(\mathcal{Y})\\
&\iff \vect(\mathcal{X} \times_3 \mathbf{A}+\mathcal{X} \times_2 \mathbf{B}+\mathcal{X} \times_1 \mathbf{C})=\vect(\mathcal{Y})\\
&\iff \mathcal{X} \times_3 \mathbf{A}+\mathcal{X} \times_2 \mathbf{B}+\mathcal{X} \times_1 \mathbf{C}=\mathcal{Y}.\\
\end{array}
\end{equation}
\section{Main Theorems}
We propose a new tensor-based formula for solving linear systems involving the Kronecker sum of matrices. This formula represents a direct solver for both the matrix and tensor forms of the Sylvester equation. Initially, the formula is derived for the case of the Kronecker sum of normal matrices and is subsequently generalized to accommodate the Kronecker sum of arbitrary matrices. In this section, we will study the case $\mathbf{\Lambda}^{(1)}_i+\mathbf{\Lambda}^{(2)}_j+\mathbf{\Lambda}^{(3)}_k \ne 0,~ \forall i,j,k$ to guarantee that the matrix $\mathbf{\Lambda}^{(3)} \oplus \mathbf{\Lambda}^{(2)} \oplus \mathbf{\Lambda}^{(1)}$ is invertible.
\subsection{Kronecker sum of normal matrices}
\begin{lemma}\label{lemma2}
If  $\left\{\mathbf{\Lambda}^{(i)}=\diag(\lambda^{(i)}_1, \dots, \lambda^{(i)}_{n_i})\right\}_{i=1}^3$ is a set of diagonal matrices, then
$$\left(\mathbf{\Lambda}^{(3)} \oplus \mathbf{\Lambda}^{(2)} \oplus \mathbf{\Lambda}^{(1)}\right)^{-1}=\diag(\vect(\cal{C}))$$ where $\cal C$ is $n_1 \times n_2 \times n_3$ tensor such that  $C_{i,j,k}=\frac{1}{\lambda^{(1)}_i+\lambda^{(2)}_j+{\lambda}^{(3)}_k}$.
\end{lemma}
\proof
Since the matrices $\left(\mathbf{\Lambda}^{(3)} \oplus \mathbf{\Lambda}^{(2)} \oplus \mathbf{\Lambda}^{(1)}\right)^{-1}$ and $\diag(\vect(\cal C))$ are diagonal, then it is enough to show that they have the same diagonal entries. But ${(\diag(\vect(\cal C)))}_{p,p}$ is the $p$-th entry of $\vect(\cal C)$ that maps to the entry $C_{i,j,k}$ such that $p=i+(j-1)n_1+(k-1)n_1 n_2$.

On the other hand,  $\left(\mathbf{\Lambda}^{(3)} \oplus \mathbf{\Lambda}^{(2)} \oplus \mathbf{\Lambda}^{(1)}\right)^{-1}_{p,p}=\frac{1}{\left(\mathbf{\Lambda}^{(3)} \oplus \mathbf{\Lambda}^{(2)} \oplus \mathbf{\Lambda}^{(1)}\right)_{p,p}}$ and
\begin{equation} 
\begin{array}{rl}
\left(\mathbf{\Lambda}^{(3)} \oplus \mathbf{\Lambda}^{(2)} \oplus \mathbf{\Lambda}^{(1)}\right)_{p,p}&=\left(\mathbf{\Lambda}^{(3)} \oplus \mathbf{\Lambda}^{(2)} \oplus \mathbf{\Lambda}^{(1)}\right)_{i+(j-1)n_1+(k-1)n_1 n_2,i+(j-1)n_1+(k-1)n_1 n_2}\\
&=\left(\mathbf{\Lambda}^{(3)} \otimes \mathbf{I}_{n_2} \otimes \mathbf{I}_{n_1}\right)_{i+(j-1)n_1+(k-1)n_1 n_2,i+(j-1)n_1+(k-1)n_1 n_2}\\
&+\left(\mathbf{I}_{n_3}\otimes \mathbf{\Lambda}^{(2)}  \otimes \mathbf{I}_{n_1}\right)_{i+(j-1)n_1+(k-1)n_1 n_2,i+(j-1)n_1+(k-1)n_1 n_2}\\
&+\left(\mathbf{I}_{n_3} \otimes \mathbf{I}_{n_2} \otimes \mathbf{\Lambda}^{(1)} \right)_{i+(j-1)n_1+(k-1)n_1 n_2,i+(j-1)n_1+(k-1)n_1 n_2}\\
&=\mathbf{\Lambda}^{(3)}_{k,k}+\mathbf{\Lambda}^{(2)}_{j,j}+\mathbf{\Lambda}^{(1)}_{i,i}\\
&=\mathbf{\lambda}^{(1)}_i+\mathbf{\lambda}^{(2)}_j+\mathbf{\lambda}^{(3)}_k.
\end{array}
\end{equation}
\endproof
\begin{theorem}\label{th1}
Let $\mathbf{A}^{(i)}$ be  normal matrices  in $\mathbb{C}^{n_i \times n_i}$, then the solution of the linear system involving Kronecker sum
\begin{equation}\label{KSLS}
\left(\mathbf{A}^{(3)} \oplus \mathbf{A}^{(2)} \oplus \mathbf{A}^{(1)}\right)\mathbf{x}=\mathbf{y},
\end{equation}
can be given as
\begin{equation}
\begin{array}{c}
\mathcal{X}=\left(\mathcal{C} \odot (\mathcal{Y} \times_1 \mathbf{U}^{(1)^*} \times_2 \mathbf{U}^{(2)^*} \times_3 \mathbf{U}^{(3)^*}\right)\times_1 \mathbf{U}^{(1)} \times_2 \mathbf{U}^{(2)} \times_3 \mathbf{U}^{(3)},\\
\end{array}
\end{equation}
where $\mathbf{U}^{(i)} \mathbf{\Lambda}^{(i)} \mathbf{U}^{(i)^*}$ be the eigendecomposition of $\mathbf{A}^{(i)}$ with unitary matrices $\mathbf{U}^{(i)}$  for $i=1, 2, 3$ and $\mathcal{X}, \mathcal{Y}$ and $\mathcal{C}$ are the $n_1 \times n_2 \times n_3$ tensors, such that $\mathbf{x}=\vect(\mathcal{X})$, $\mathbf{y}=\vect(\mathcal{Y})$, 
$C_{i,j,k}=\frac{1}{\lambda^{(1)}_i+\lambda^{(2)}_j+\lambda^{(3)}_k}$
 and 
$\mathbf{\Lambda}^{(i)}=\diag(\lambda^{(i)}_1, \dots, \lambda^{(i)}_{n_i}),$  for $i=1, 2, 3$.
\end{theorem}
\proof
Consider the linear system (\ref{KSLS}) where
\begin{equation}
\mathbf{A}^{(3)} \oplus \mathbf{A}^{(2)} \oplus \mathbf{A}^{(1)}=\mathbf{A}^{(3)}\otimes \mathbf{I}_{n_2}\otimes \mathbf{I}_{n_1}+\mathbf{I}_{ n_3}\otimes \mathbf{A}^{(2)}\otimes \mathbf{I}_{n_1}+\mathbf{I}_{n_3}\otimes \mathbf{I}_{n_2}\otimes \mathbf{A}^{(1)}.\\
\end{equation}
Using the eigendecomposition of $\mathbf{A}^{(i)}$, the linear system (\ref{KSLS}) can be written as
\begin{equation}
\left(\left(\mathbf{U}^{(3)}\otimes \mathbf{U}^{(2)} \otimes \mathbf{U}^{(1)}\right) \left(\mathbf{\Lambda}^{(3)} \oplus \mathbf{\Lambda}^{(2)} \oplus \mathbf{\Lambda}^{(1)}\right) \left(\mathbf{U}^{(3)}\otimes \mathbf{U}^{(2)} \otimes \mathbf{U}^{(1)}\right)^*\right) \mathbf{x}=\mathbf{y}
\end{equation}
Then
\begin{equation}\label{MEQ1}
\left(\mathbf{U}^{(3)}\otimes \mathbf{U}^{(2)} \otimes \mathbf{U}^{(1)}\right)^* \mathbf{x}=\left(\mathbf{\Lambda}^{(3)} \oplus \mathbf{\Lambda}^{(2)} \oplus \mathbf{\Lambda}^{(1)}\right)^{-1} \left(\mathbf{U}^{(3)}\otimes \mathbf{U}^{(2)} \otimes \mathbf{U}^{(1)}\right)^* \mathbf{y} \\
\end{equation}
If we define the $n_1 \times n_2 \times n_3$ tensors $\mathcal{X}$ and $\mathcal{Y}$ such that $\mathbf{x}=\vect(\mathcal{X})$ and $\mathbf{y}=\vect(\mathcal{Y})$, $\hat{\mathcal{Y}}=\mathcal{Y} \times_1 \mathbf{U}^{(1)^*} \times_2 \mathbf{U}^{(2)^*} \times_3 \mathbf{U}^{(3)^*}$ and $\hat{\mathcal{X}}=\mathcal{X} \times_1 \mathbf{U}^{(1)^*} \times_2 \mathbf{U}^{(2)^*} \times_3 \mathbf{U}^{(3)^*}$, then the above formula can be simplified further to the simple version as follows:
\begin{equation}
\begin{array}{ll}
\vect(\hat{\mathcal{X}})&=\left(\mathbf{U}^{(3)}\otimes \mathbf{U}^{(2)} \otimes 
\mathbf{U}^{1}\right)^* \vect(\mathcal{X})\\
\vect(\mathcal{C}\odot \hat{\mathcal{Y}})&=\vect(\mathcal{C}) \odot \vect(\hat{\mathcal{Y}})\\
&=\diag(\vect(\mathcal{C})) \vect(\hat{\mathcal{Y}})\\
&=\left(\mathbf{\Lambda}^{(3)} \oplus \mathbf{\Lambda}^{(2)} \oplus \mathbf{\Lambda}^{(1)}\right)^{-1}\left(\mathbf{U}^{(3)}\otimes \mathbf{U}^{(2)} \otimes \mathbf{U}^{(1)}\right)^* \vect({\mathcal{Y}})\\
\end{array}
\end{equation}
Then using (\ref{MEQ1}), we have $\hat{\mathcal{X}}=\mathcal{C} \odot \hat{\mathcal{Y}}.$
Finally, since  ${\mathcal{X}}=\hat{\mathcal{X}} \times_1 \mathbf{U}^{(1)} \times_2 \mathbf{U}^{(2)} \times_3 \mathbf{U}^{(3)},$ then 
\begin{equation}
\begin{array}{rl}
{\mathcal{X}}&=\left(\mathcal{C} \odot \hat{\mathcal{Y}}\right)\times_1 \mathbf{U}^{(1)} \times_2 \mathbf{U}^{(2)} \times_3 \mathbf{U}^{(3)},\\
&=\left(\mathcal{C} \odot ({\mathcal{Y}} \times_1 \mathbf{U}^{(1)^*} \times_2 \mathbf{U}^{(2)^*} \times_3 \mathbf{U}^{(3)^*})\right)\times_1 \mathbf{U}^{(1)} \times_2 \mathbf{U}^{(2)} \times_3 \mathbf{U}^{(3)}.\\
\end{array}
\end{equation}
\endproof
\begin{corollary}\label{co1}
Let $\mathbf{A}^{(i)}$ be  normal matrices  in $\mathbb{C}^{n_i \times n_i}$, then the solution of the linear system involving Kronecker sum
\begin{equation}
\left(\mathbf{A}^{(2)} \oplus \mathbf{A}^{(1)}\right)\mathbf{x}=\mathbf{y},
\end{equation}
can be given as
\begin{equation}
\begin{array}{c}
\mathbf{X}=\mathbf{U}^{(1)} \left(\mathbf{C} \odot (\mathbf{U}^{(1)^*} \mathbf{Y} \mathbf{U}^{(2)})\right){\mathbf{U}^{(2)}}^*\\
\end{array}
\end{equation}
where $\mathbf{U}^{(i)} \mathbf{\Lambda}^{(i)} \mathbf{U}^{(i)^*}$ be the eigendecomposition of $\mathbf{A}^{(i)}$ with unitary matrices $\mathbf{U}^{(i)}$  for $i=1, 2$ and $\mathbf{X}, \mathbf{Y}$ and $\mathbf{C}$ are the $n_1 \times n_2$ matrices, such that $\mathbf{x}=\vect(\mathbf{X})$, $\mathbf{y}=\vect(\mathbf{Y}),$ $C_{ij}=\frac{1}{\lambda^{(1)}_i+\lambda^{(2)}_j}$ and $\mathbf{\Lambda}^{(i)}=\diag(\lambda^{(i)}_1, \dots, \lambda^{(i)}_{n_i}),$  for $i=1, 2$.
\end{corollary}
\proof
Use Theorem  (\ref{th1}) with the facts $\mathbf{Y} \times_1 \mathbf{B}=\mathbf{B}\mathbf{Y}$ and $\mathbf{Y} \times_2 \mathbf{C}=\mathbf{Y} \mathbf{C}^T$ for any  in $\mathbf{B} \in \mathbb{C}^{n_1 \times n_1}$ and $\mathbf{C} \in \mathbb{C}^{n_2 \times n_2}$.
\endproof
\subsection{Kronecker sum of matrices}
\begin{lemma}\label{lemma3}
If  $\left\{\mathbf{T}^{(i)}\right\}_{i=1}^3$ is a set of $n_i \times n_i$ upper triangular matrices, then
$$\left(\mathbf{T}^{(3)} \oplus \mathbf{T}^{(2)} \oplus \mathbf{T}^{(1)}\right)^{-1}=\left(\sum\limits_{j=0}^{d-1}{(-1)}^j \mathbf{N}^j \right)\left(\mathbf{\Lambda}^{(3)} \oplus \mathbf{\Lambda}^{(2)} \oplus \mathbf{\Lambda}^{(1)}\right)^{-1},$$ 
where 
$\mathbf{\Lambda}^{(i)}=\diag(\mathbf{T}^{(i)})$ is a diagonal matrix, $\mathbf{N}=\left(\mathbf{\Lambda}^{(3)} \oplus \mathbf{\Lambda}^{(2)} \oplus \mathbf{\Lambda}^{(1)}\right)^{-1}\left(\tilde{ \mathbf{T}}^{(3)} \oplus \tilde{\mathbf{T}}^{(2)} \oplus \tilde{\mathbf{T}}^{(1)}\right)$  is a nilpotent matrix with index $d \leq n_1n_2n_3$ and $\tilde{\mathbf{T}}^{(i)}=\mathbf{T}^{(i)}-\mathbf{\Lambda}^{(i)}$ for $i=1, 2, 3.$
\end{lemma}
\proof
Let $\mathbf{\Lambda}^{(i)}$ be the diagonal matrix formed from the diagonal of $\mathbf{T}^{(i)}$ and 
$\tilde{\mathbf{T}}^{(i)}$ be the strictly upper triangular matrix formed by setting the diagonal of $\mathbf{T}^{(i)}$ to zero for $i=1, 2, 3.$ Then $\mathbf{T}^{(i)}$ can be written as
$$\mathbf{T}^{(i)}=\mathbf{\Lambda}^{(i)}+ \tilde{ \mathbf{T}}^{(i)},~i=1,2,3.$$

Consequently, $\left(\mathbf{T}^{(3)} \oplus \mathbf{T}^{(2)} \oplus \mathbf{T}^{(1)}\right)$ can be written as
\begin{equation}
\begin{array}{rl}
\left(\mathbf{T}^{(3)} \oplus \mathbf{T}^{(2)} \oplus \mathbf{T}^{(1)}\right)&
=\left(\mathbf{T}^{(3)} \otimes \mathbf{I}_{n_2} \otimes \mathbf{I}_{n_1} \right)
+\left(\mathbf{I}_{n_3} \otimes \mathbf{T}^{(2)} \otimes \mathbf{I}_{n_1} \right)
+\left(\mathbf{I}_{n_3} \otimes \mathbf{I}_{n_2} \otimes \mathbf{T}^{(1)}\right)\\
&=\left(\mathbf{\Lambda}^{(3)} \otimes \mathbf{I}_{n_2} \otimes \mathbf{I}_{n_1} \right)
+\left(\mathbf{I}_{n_3} \otimes \mathbf{\Lambda}^{(2)} \otimes \mathbf{I}_{n_1} \right)
+\left(\mathbf{I}_{n_3} \otimes \mathbf{I}_{n_2} \otimes \mathbf{\Lambda}^{(1)}\right)\\
&+\left(\tilde{\mathbf{T}}^{(3)} \otimes \mathbf{I}_{n_2} \otimes \mathbf{I}_{n_1} \right)
+\left(\mathbf{I}_{n_3} \otimes \tilde{\mathbf{T}}^{(2)} \otimes \mathbf{I}_{n_1} \right)
+\left(\mathbf{I}_{n_3} \otimes \mathbf{I}_{n_2} \otimes \tilde{\mathbf{T}}^{(1)}\right)\\
&=\left(\mathbf{\Lambda}^{(3)} \oplus \mathbf{\Lambda}^{(2)} \oplus \mathbf{\Lambda}^{(1)}\right)
+\left(\tilde{\mathbf{T}}^{(3)} \oplus \tilde{ \mathbf{T}}^{(2)} \oplus \tilde{ \mathbf{T}}^{(1)}\right)\\
&=\left(\mathbf{\Lambda}^{(3)} \oplus \mathbf{\Lambda}^{(2)} \oplus \mathbf{\Lambda}^{(1)}\right)
\left(\mathbf{I}_{n_1 n_2 n_3}+\mathbf{N} \right).\\
\end{array}
\end{equation}
Therefore, the inverse of $\left(\mathbf{T}^{(3)} \oplus \mathbf{T}^{(2)} \oplus \mathbf{T}^{(1)}\right)$ can be given as
\begin{equation}
\begin{array}{rl}
\left(\mathbf{T}^{(3)} \oplus \mathbf{T}^{(2)} \oplus \mathbf{T}^{(1)}\right)^{-1}
=\left(\mathbf{I}_{n_1 n_2 n_3}+\mathbf{N} \right)^{-1}
\left(\mathbf{\Lambda}^{(3)} \oplus \mathbf{\Lambda}^{(2)} \oplus \mathbf{\Lambda}^{(1)}\right)^{-1}.\\
\end{array}
\end{equation}
Finally, since $\mathbf{N}$ is a nilpotent matrix with index $d \leq n_1n_2n_3$, then
\begin{equation}
\begin{array}{rl}
\left(\mathbf{I}_{n_1 n_2 n_3}+\mathbf{N} \right)^{-1}=\mathbf{I}_{n_1 n_2 n_3}+\sum\limits_{j=1}^{d-1}{(-1)}^j \mathbf{N}^j.
\end{array}
\end{equation}
\endproof

\begin{theorem}\label{th2}
Let $\mathbf{A}^{(i)}$ be matrices  in $\mathbb{C}^{n_i \times n_i}$, then the solution of the linear system involving Kronecker sum
\begin{equation}\label{KSLS2}
\left(\mathbf{A}^{(3)} \oplus \mathbf{A}^{(2)} \oplus \mathbf{A}^{(1)}\right)\mathbf{x}=\mathbf{y},
\end{equation}
can be given as
\begin{equation}
\begin{array}{c}
\mathcal{X}=\left(\sum\limits_{j=0}^{n_1n_2n_3-1}{(-1)}^j \mathcal{K}_j\right)\times_1 \mathbf{U}^{(1)} \times_2 \mathbf{U}^{(2)} \times_3 \mathbf{U}^{(3)},\\
\end{array}
\end{equation}
such that the tensors $\mathcal{K}_{j}$ can be evaluated using the recursive formula
\begin{equation}
\begin{array}{l}
\mathcal{K}_0=\mathcal{C} \odot  \left(\mathcal{Y} \times_1 \mathbf{U}^{(1)^*} \times_2 \mathbf{U}^{(2)^*} \times_3 \mathbf{U}^{(3)^*}\right),\\
\mathcal{K}_{j+1}=\mathcal{C} \odot \left(\mathcal{K}_j \times_1 \tilde{\mathbf{T}}^{(1)}+\mathcal{K}_j \times_2 \tilde{\mathbf{T}}^{(2)} +\mathcal{K}_j \times_3 \tilde{\mathbf{T}}^{(3)}\right),\\
\end{array}
\end{equation}
where $\mathbf{U}^{(i)} \mathbf{T}^{(i)} \mathbf{U}^{(i)^*}$ be the Schur decomposition of $\mathbf{A}^{(i)}$ with unitary matrices $\mathbf{U}^{(i)}$  for $i=1,2, 3$ and $\mathcal{X}, \mathcal{Y}$ and $\mathcal{C}$ are the $n_1 \times n_2 \times n_3$ tensors, such that $\mathbf{x}=\vect(\mathcal{X})$, $\mathbf{y}=\vect(\mathcal{Y}),$ $C_{i,j,k}=\frac{1}{\lambda^{(1)}_i+\lambda^{(2)}_j+\lambda^{(3)}_k}$ and $\mathbf{\Lambda}^{(i)}=\diag(\mathbf{T}^{(i)})=\diag(\lambda^{(i)}_1, \dots, \lambda^{(i)}_{n_i})$  and $\tilde{\mathbf{T}}^{(i)}=\mathbf{T}^{(i)}-\mathbf{\Lambda}^{(i)}$ for $i=1, 2, 3.$
\end{theorem}
\proof
Consider the linear system (\ref{KSLS2}) where
\begin{equation}
\mathbf{A}^{(3)} \oplus \mathbf{A}^{(2)} \oplus \mathbf{A}^{(1)}=\mathbf{A}^{(3)}\otimes \mathbf{I}_{n_2}\otimes \mathbf{I}_{n_1}+\mathbf{I}_{ n_3}\otimes \mathbf{A}^{(2)}\otimes \mathbf{I}_{n_1}+\mathbf{I}_{n_3}\otimes \mathbf{I}_{n_2}\otimes \mathbf{A}^{(1)}.\\
\end{equation}
Using the Schur decomposition of $\mathbf{A}^{(i)}$,  the linear system (\ref{KSLS2}) can be written as
\begin{equation}
\left(\left(\mathbf{U}^{(3)}\otimes \mathbf{U}^{(2)} \otimes \mathbf{U}^{(1)}\right) \left(\mathbf{T}^{(3)} \oplus \mathbf{T}^{(2)} \oplus \mathbf{T}^{(1)}\right) \left(\mathbf{U}^{(3)}\otimes \mathbf{U}^{(2)} \otimes \mathbf{U}^{(1)}\right)^*\right) \mathbf{x}=\mathbf{y}
\end{equation}
Using Lemma (\ref{lemma3}), then
\begin{equation}\label{MEQ2}
\begin{array}{ll}
\left(\mathbf{U}^{(3)}\otimes \mathbf{U}^{(2)} \otimes \mathbf{U}^{(1)}\right)^* \mathbf{x}&=\left(\mathbf{T}^{(3)} \oplus \mathbf{T}^{(2)} \oplus \mathbf{T}^{(1)}\right)^{-1} \left(\mathbf{U}^{(3)}\otimes \mathbf{U}^{(2)} \otimes \mathbf{U}^{(1)}\right)^*  \mathbf{y}\\
&=\left(\sum\limits_{j=0}^{n_1n_2n_3-1}{(-1)}^j \mathbf{N}^j \right)\left(\mathbf{\Lambda}^{(3)} \oplus \mathbf{\Lambda}^{(2)} \oplus \mathbf{\Lambda}^{(1)}\right)^{-1}\left(\mathbf{U}^{(3)}\otimes \mathbf{U}^{(2)} \otimes \mathbf{U}^{(1)}\right)^*  \mathbf{y},\\
\end{array}
\end{equation}
where $\mathbf{N}=\left(\mathbf{\Lambda}^{(3)} \oplus \mathbf{\Lambda}^{(2)} \oplus \mathbf{\Lambda}^{(1)}\right)^{-1}\left(\tilde{\mathbf{T}}^{(3)} \oplus \tilde{\mathbf{T}}^{(2)} \oplus \tilde{\mathbf{T}}^{(1)}\right).$

If we define the $n_1 \times n_2 \times n_3$ tensors $\mathcal{X}$ and $\mathcal{Y}$ such that $\mathbf{x}=\vect(\mathcal{X})$ and $\mathbf{y}=\vect(\mathcal{Y})$, $\hat{\mathcal{Y}}=\mathcal{Y} \times_1 \mathbf{U}^{(1)^*} \times_2 \mathbf{U}^{(2)^*} \times_3 \mathbf{U}^{(3)^*}$ and $\hat{\mathcal{X}}=\mathcal{X} \times_1 \mathbf{U}^{(1)^*} \times_2 \mathbf{U}^{(2)^*} \times_3 \mathbf{U}^{(3)^*}$, then the above formula can be simplified further using Lemma (\ref{lemma1}) and Lemma (\ref{lemma2}) to the simple version as follows:
\begin{equation}
\begin{array}{ll}
\vect(\hat{\mathcal{X}})&=\left(\mathbf{U}^{(3)}\otimes \mathbf{U}^{(2)} \otimes \mathbf{U}^{1}\right)^* \vect({\mathcal{X}})\\
\vect(\mathcal{K}_0)&=\vect(\mathcal{C} \odot \hat{\mathcal{Y}})=\vect(\mathcal{C}) \odot \vect(\hat{\mathcal{Y}})\\
&=\diag(\vect(\mathcal{C})) \vect(\hat{\mathcal{Y}})\\
&=\left(\mathbf{\Lambda}^{(3)} \oplus \mathbf{\Lambda}^{(2)} \oplus \mathbf{\Lambda}^{(1)}\right)^{-1}\left(\mathbf{U}^{(3)}\otimes \mathbf{U}^{(2)} \otimes \mathbf{U}^{(1)}\right)^* \vect(\mathcal{Y})\\
\vect(\mathcal{K}_{j+1})&=\vect(\mathcal{C} \odot \left(\mathcal{K}_j \times_1 \tilde{\mathbf{T}}^{(1)}+\mathcal{K}_j \times_2 \tilde{\mathbf{T}}^{(2)}+\mathcal{K}_j \times_3 \tilde{\mathbf{T}}^{(3)}\right)\\
&=\vect(\mathcal{C}) \odot vec\left(\mathcal{K}_j \times_1 \tilde{\mathbf{T}}^{(1)}+\mathcal{K}_j \times_2 \tilde{\mathbf{T}}^{(2)}+\mathcal{K}_j \times_3 \tilde{\mathbf{T}}^{(3)}\right)\\
&=\diag(\vect(\mathcal{C}))vec\left(\mathcal{K}_j \times_1 \tilde{\mathbf{T}}^{(1)}+\mathcal{K}_j \times_2 \tilde{\mathbf{T}}^{(2)}+\mathcal{K}_j \times_3 \tilde{\mathbf{T}}^{(3)}\right)\\
&=\diag(\vect(\mathcal{C})) \left(\mathbf{I}_{n_3}\otimes \mathbf{I}_{n_2}\otimes \tilde{\mathbf{T}}^{(1)}+\mathbf{I}_{ n_3}\otimes \tilde{\mathbf{T}}^{(2)}\otimes \mathbf{I}_{n_1}+\tilde{\mathbf{T}}^{(3)}\otimes \mathbf{I}_{n_2}\otimes \mathbf{I}_{n_1}\right) \vect(\mathcal{K}_j)\\
&=\diag(\vect(\mathcal{C})) \left(\tilde{\mathbf{T}}^{(3)}\oplus\tilde{\mathbf{T}}^{(2)}\oplus\tilde{\mathbf{T}}^{(1)}\right) \vect(\mathcal{K}_j)\\
&=\left(\mathbf{\Lambda}^{(3)} \oplus \mathbf{\Lambda}^{(2)} \oplus \mathbf{\Lambda}^{(1)}\right)^{-1} \left(\tilde{\mathbf{T}}^{(3)}\oplus \tilde{\mathbf{T}}^{(2)}\oplus \tilde{\mathbf{T}}^{(1)}\right) \vect(\mathcal{K}_j)\\
&=\mathbf{N} \vect(\mathcal{K}_j)\\
&=\mathbf{N}^{j+1} \vect(\mathcal{K}_0)\\
\end{array}
\end{equation}
Therefore
\begin{equation}
\begin{array}{ll}
\vect\left(\sum\limits_{j=0}^{n_1n_2n_3-1}{(-1)}^j \mathcal{K}_j \right)&=\sum\limits_{j=0}^{n_1n_2n_3-1}{(-1)}^j \vect(\mathcal{K}_j)\\
&=\sum\limits_{j=0}^{n_1n_2n_3-1}{(-\mathbf{N})}^j \vect(\mathcal{K}_0)\\
&=\sum\limits_{j=0}^{n_1n_2n_3-1}{(-\mathbf{N})}^j \left(\mathbf{\Lambda}^{(3)} \oplus \mathbf{\Lambda}^{(2)} \oplus \mathbf{\Lambda}^{(1)}\right)^{-1}\left(\mathbf{U}^{(3)}\otimes \mathbf{U}^{(2)} \otimes \mathbf{U}^{(1)}\right)^* \vect({\mathcal{Y}})\\
\end{array}
\end{equation}
Then using (\ref{MEQ2}), we have $$\hat{\mathcal{X}}=\sum\limits_{j=0}^{n_1n_2n_3-1}{(-1)}^j \mathcal{K}_j.$$
Finally, since  $\mathcal{X}=\hat{\mathcal{X}} \times_1 \mathbf{U}^{(1)} \times_2 \mathbf{U}^{(2)} \times_3 \mathbf{U}^{(3)},$ then 
\begin{equation}
\begin{array}{rl}
\mathcal{X}&=\left(\sum\limits_{j=0}^{n_1n_2n_3-1}{(-1)}^j \mathcal{K}_j \right)\times_1 \mathbf{U}^{(1)} \times_2 \mathbf{U}^{(2)} \times_3 \mathbf{U}^{(3)}.\\
\end{array}
\end{equation}
\endproof
\begin{remark}
The formula in Theorem  (\ref{th2})  can be generalized to a sum of a finite number of matrices using a similar proof. Note that a matrix has a spectral radius less than one if and only if the successive powers of a matrix converge to the zero matrix. Therefore, since $\mathbf{N}$ is a nilpotent matrix and $\vect(\mathcal{K}_{j+1})=\mathbf{N}^{j+1} \vect(\mathcal{K}_0)$, then the tensors $\mathcal{K}_{j}$ converges to zero tensor as $j$ increases, so one can approximate the answer by truncating the summation $\sum\limits_{j=0}^{n_1n_2n_3-1}{(-1)}^j \mathcal{K}_j$.
\end{remark}
\begin{corollary}\label{co2}

Let $\mathbf{A}^{(i)}$ be matrices  in $\mathbb{C}^{n_i \times n_i}$, then the solution of the linear system involving Kronecker sum
\begin{equation}\label{KSLS2}
\left( \mathbf{A}^{(2)} \oplus \mathbf{A}^{(1)}\right)\mathbf{x}=\mathbf{y},
\end{equation}
can be given as
\begin{equation}
\begin{array}{c}
\mathbf{X}=\mathbf{U}^{(1)} \left(\sum\limits_{j=0}^{n_1n_2-1}{(-1)}^j \mathcal{K}_j\right) \mathbf{U}^{(2)^*},\\
\end{array}
\end{equation}
such that the tensors $\mathcal{K}_{j}$ can be evaluated using the recursive formula
\begin{equation}
\begin{array}{l}
\mathcal{K}_0=\mathbf{C} \odot  \left( \mathbf{U}^{(1)^*}\mathbf{Y}  \mathbf{U}^{(2)}\right),\\
\mathcal{K}_{j+1}=\mathbf{C} \odot \left( \mathbf{\tilde{T}}^{(1)} \mathcal{K}_j +\mathcal{K}_j \mathbf{\tilde{T}}^{(1)^T}\right),\\
\end{array}
\end{equation}
where $\mathbf{U}^{(i)} \mathbf{T}^{(i)} \mathbf{U}^{(i)^*}$ be the Schur decomposition of $\mathbf{A}^{(i)}$ with unitary matrices $\mathbf{U}^{(i)}$  for $i=1, 2$ and $\mathbf{X}, \mathbf{Y}$ and $\mathbf{C}$ are the $n_1 \times n_2$ matrices, such that $\mathbf{x}=\vect(\mathbf{X})$, $\mathbf{y}=\vect(\mathbf{Y}),$ $C_{ij}=\frac{1}{\lambda^{(1)}_i+\lambda^{(2)}_j}$ and $\mathbf{\Lambda}^{(i)}=\diag(\mathbf{T}^{(i)})=\diag(\lambda^{(i)}_1, \dots, \lambda^{(i)}_{n_i})$  and $\tilde{\mathbf{T}}^{(i)}=\mathbf{T}^{(i)}-\mathbf{\Lambda}^{(i)}$ for $i=1, 2.$
\end{corollary}
\proof
Use Theorem  (\ref{th2}) with the facts $\mathbf{Y} \times_1 \mathbf{B}=\mathbf{B} \mathbf{Y}$ and $\mathbf{Y} \times_2 \mathbf{C}=\mathbf{Y} \mathbf{C}^T$ for any  in $\mathbf{B} \in \mathbb{C}^{n_1 \times n_1}$ and $\mathbf{C} \in \mathbb{C}^{n_2 \times n_2}$.
\endproof
\section{Application}
\subsection{Explicit solution of 2D discrete Poisson equation}
The 2D Poisson equation defined on a domain $\Omega=[0,1]\times [0,1]$
with Dirichlet boundary conditions is given by
\begin{equation}
u_{xx}+u_{yy}=f, u(x,y)|_{\partial \Omega}=g,
\end{equation}
which can be discretized using finite differences to obtain a system of linear equations of the form
\begin{equation}\label{2DPoisson}
\mathbf{A}_2 \mathbf{u}=\mathbf{b}.
\end{equation}
where $\mathbf{A}_2$ is a coefficient matrix, $\mathbf{u}$ the discrete solution and $\mathbf{b}$
contains terms from the source $f$ and the boundary conditions $g$.

The Kronecker sum representation of discretized Laplacian in 2D is
$$\mathbf{A}_2=\mathbf{A}_1\oplus \mathbf{A}_1=\mathbf{I}_N\otimes \mathbf{A}_1+\mathbf{A}_1\otimes \mathbf{I}_N,$$ where
$\mathbf{A}_1$ is the one-dimensional second-derivative approximation
\begin{equation}\label{T}
\mathbf{A}_1= \begin{pmatrix}
2& -1 & 0 & \cdots & 0 \\
-1  & 2 &-1 & \ddots & \vdots \\
0 & -1 & 2 & \ddots & 0\\
\vdots & \ddots & \ddots & \ddots & -1 \\
0 & \cdots & 0 & -1  & 2
\end{pmatrix}_{N\times N}.
\end{equation}
Since $\mathbf{A}_1$ is Unitarily diagonalizable, then $\mathbf{A}_1=\mathbf{U}\mathbf{D}\mathbf{U}^*$ for some diagonal matrix $\mathbf{D}=\diag(\lambda_1,\dots,\lambda_N)$ and unitary matrix $\mathbf{U}=\left[\mathbf{s}_1 \mathbf{s}_2 \cdots \mathbf{s}_N\right]_{N\times N}$ where
\begin{equation}\label{lambda}
\lambda_i=4 \sin^2\left( \frac{i \pi}{2(N+1)}\right),
\end{equation}
and 
\begin{equation}\label{si}
\mathbf{s}_i=\sqrt{\frac{2}{N+1}}\left[\sin\left( \frac{i \pi}{N+1}\right)~\sin\left( \frac{2i \pi}{N+1}\right) ~\dots~ \sin\left( \frac{N i \pi}{N+1}\right)\right]^T.
\end{equation}
If we define the $N \times N$ matrices $\mathbf{X}$ and $\mathbf{Y}$ such that $\mathbf{u}=\vect(\mathbf{X})$ and $\mathbf{b}=\vect(\mathbf{Y})$, then the solution can be given explicitly using Corollary (\ref{co1}) as
\begin{equation}
\begin{array}{cc}
\mathbf{X}=\mathbf{U} \left(\mathbf{C} \odot (\mathbf{U}^* \mathbf{Y} \mathbf{U})\right)\mathbf{U}^*,& C_{i,j}=\frac{1}{\lambda_i+\lambda_j}.\\
\end{array}
\end{equation}
So, one can solve the system (\ref{2DPoisson}) as follows:
\begin{enumerate}
\item[Step 1:]  Compute the diagonal matrix $\mathbf{D}=\diag(\lambda_1,\dots,\lambda_N)$ and unitary matrix $\mathbf{U}=\left[\mathbf{s}_1 \mathbf{s}_2 \cdots \mathbf{s}_N\right]_{N\times N}$ using (\ref{lambda}) and (\ref{si}).
\item[Step 2:]  Compute $N \times N$ matrices $\mathbf{C}$ and $\mathbf{Y}$ such that $C_{i,j}=\frac{1}{\lambda_i+\lambda_j}$ and $\mathbf{b}=\vect(\mathbf{Y}).$
\item[Step 3:] Compute the matrix $\mathbf{T}_1=\mathbf{U}^* \mathbf{Y} \mathbf{U}.$
\item[Step 4:] Compute the matrix  $\mathbf{T}_2=\mathbf{C} \odot \mathbf{T}_1$.
\item[Step 5:] Compute the matrix $\mathbf{X}=\mathbf{U} \mathbf{T}_2 \mathbf{U}^*.$
\item[Step 6:] Compute the solution $\mathbf{u}=\vect(\mathbf{X})$.
\end{enumerate}
The computational cost of the algorithm based on Table \ref{table:1} is $\mathcal {O}(N^3)=\mathcal {O}(n^{3/2})$. It should be mentioned here that the Fast Poisson Solver based on diagonalization given in \cite{Lyche2020} has the same cost $\mathcal {O}(n^{3/2})$.
\begin{table}[H]
\centering
\begin{tabular}{|c| c| c| c|} 
 \hline
 Step &  $+, -$ & $\times,  /$ & Total \\
\hline\hline
2 & $N^2$ & $N^2$ & $2N^2$ \\ 
3 &$2N^2(N-1)$  & $2N^3$ & $4N^3-2N^2$ \\ 
4 & 0 & $N^2$ & $N^2$  \\ 
5 &$2N^2(N-1)$  & $2N^3$ & $4N^3-2N^2$ \\ 
  \hline
 & & & $8N^3-N^2$  \\ 
 \hline
\end{tabular}
\caption{The computational cost of the algorithm for 2D Poisson equation}
\label{table:1}
\end{table}
\begin{remark}
We can compute $\mathbf{U}$ and $\mathbf{D}$ in step 1 without using loops. Using outer products and raising a matrix element by element to a power as follows:
\begin{equation}
\begin{array}{l}
(\lambda_1,\dots,\lambda_N)=4 \sin\left( \frac{\pi}{2(N+1)}
\begin{bmatrix}
1&2 &\cdots& N
\end{bmatrix}\right)^2,\\
\mathbf{U}=\sqrt{\frac{2}{N+1}}\sin\left( \frac{\pi}{N+1}
\begin{bmatrix}
1 \\
2 \\
\vdots \\
N
\end{bmatrix}
\begin{bmatrix}
1&2 &\cdots& N
\end{bmatrix}
\right).\\
\end{array}
\end{equation}
\end{remark}
\begin{figure}[H]
\centering
\includegraphics[totalheight=10cm]{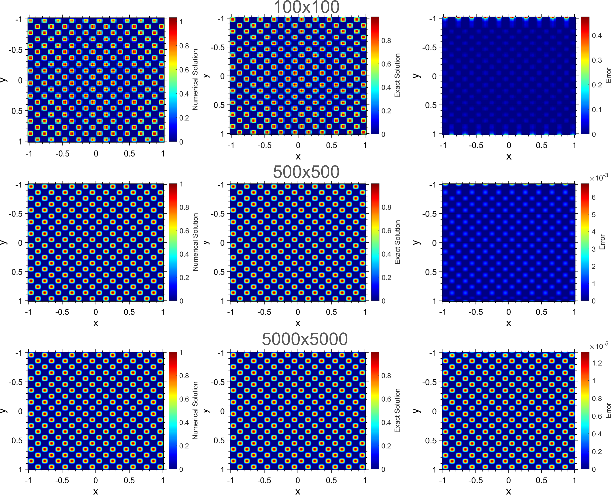}
\caption{Comparison of the numerical and exact solutions at different resolutions for the $2D$ Poisson equation on the domain  $[-1,1]\times [-1,1]$ with zero Dirichlet boundary conditions . The load term is $-200 \pi^2 \sin(10\pi x) \sin(10\pi y)$. Results are shown for different grid resolutions.}
\label{Figure 2}
\end{figure}
\subsection{Explicit solution of 3D discrete Poisson equation}
Similarly, the 3D Poisson equation defined on a domain $\Omega=[0,1]\times [0,1] \times [0,1]$ with Dirichlet boundary conditions
\begin{equation}
u_{xx}+u_{yy}+u_{zz}=f, u(x,y,z)|_{\partial \Omega}=g,
\end{equation}
can be discretized, using finite differences, to obtain a system of linear equations of the form
\begin{equation}\label{3DPoisson}
\mathbf{A}_3 \mathbf{u}=\mathbf{b}.
\end{equation}
The Kronecker sum representation of discretized Laplacian $\mathbf{A}_3$ in 3D is
$$\mathbf{A}_3=\mathbf{A}_1 \oplus \mathbf{A}_1 \oplus \mathbf{A}_1=\mathbf{I}_N\otimes \mathbf{I}_N\otimes \mathbf{A}_1+\mathbf{I}_N\otimes \mathbf{A}_1\otimes \mathbf{I}_N+\mathbf{A}_1\otimes \mathbf{I}_N\otimes \mathbf{I}_N$$ where $\mathbf{A}_1$ is given by equation (\ref{T}).

If we define the $N \times N \times N$ tensors $\mathcal{X}$ and $\mathcal{Y}$ such that $\mathbf{u}=\vect(\mathcal{X})$ and $\mathbf{b}=\vect(\mathcal{Y})$, then the solution can be given explicitly using Theorem (\ref{th1}) as
\begin{equation}
\begin{array}{cc}
\mathcal{X}=\left(\mathcal{C} \odot (\mathcal{Y} \times_1 \mathbf{U}^* \times_2 \mathbf{U}^* \times_3 \mathbf{U}^*)\right)\times_1 \mathbf{U} \times_2 \mathbf{U} \times_3 \mathbf{U},& C_{i,j,k}=\frac{1}{\lambda_i+\lambda_j+\lambda_k},\\
\end{array}
\end{equation}
where the tensor ${\mathcal{T}} \times _n \mathbf{V}^{(n)}$ is the $n$-mode product of the 3rd-order tensor ${\cal T} \in \mathbb{C}^{d_1  \times d_2  \times d_3 }$ and matrices $\mathbf{V}^{(1)}, \mathbf{V}^{(2)} \mathbf{V}^{(3)}$ in $\mathbb{C}^{k_1  \times d_1 } ,\mathbb{C}^{k_2  \times d_2 } ,\mathbb{C}^{k_3  \times d_3 },$  respectively, i.e., 
\begin{equation}\label{tensor mode product}
\begin{array}{cc}
({\cal T} \times _1 \mathbf{V}^{(1)})_{i_1 ,j_2 ,j_3} = \sum\limits_{j_1  = 1}^{d_1 } {\mathcal{T}}_{j_1 ,j_2 ,j_3 }V^{(1)}_ {i_1 ,j_1},\\
({\mathcal{T}} \times _2 \mathbf{V}^{(2)})_{ j_1, i_2 ,j_3} = \sum\limits_{j_2  = 1}^{d_2 } {\cal T}_{j_1 ,j_2 ,j_3 }V^{(2)}_ {i_2 ,j_2},\\
({\cal T} \times _3 \mathbf{V}^{(3)})_{j_1 ,j_2, i_3} = \sum\limits_{j_3  = 1}^{d_3 } {\mathcal{T}}_{j_1 ,j_2 ,j_3 }V^{(3)}_ {i_3,j_3}.\\
\end{array}
\end{equation}
So, one can solve the system (\ref{3DPoisson}) as follows:
\begin{enumerate}
\item[Step 1:]  Compute the diagonal matrix $\mathbf{D}=\diag(\lambda_1,\dots,\lambda_N)$ and unitary matrix $\mathbf{U}=[\mathbf{s}_1 \mathbf{s}_2 \cdots \mathbf{s}_N]_{N\times N}$ using (\ref{lambda}) and (\ref{si}).
\item[Step 2:]  Compute $N \times N \times N$ tensors $\mathcal{C}$ and $\mathcal{Y}$ such that $C_{ijk}=\frac{1}{\lambda_i+\lambda_j+\lambda_k}$ and $\mathbf{b}=\vect(\mathcal{Y})$ .
\item[Step 3:] Compute the tensor $\mathcal{T}_1=\mathcal{Y} \bigtimes_1 \mathbf{U}^* \bigtimes_2 \mathbf{U}^* \bigtimes_3 \mathbf{U}^*$ using (\ref{tensor mode product}).
\item[Step 4:] Compute the tensor $\mathcal{T}_2=\mathcal{C} \odot \mathcal{T}_1$.
\item[Step 5:] Compute the tensor $\mathcal{X}=\mathcal{T}_2 \bigtimes_1 \mathbf{U} \bigtimes_2 \mathbf{U} \bigtimes_3 \mathbf{U}$ using (\ref{tensor mode product}).
\item[Step 6:] Compute the solution $\mathbf{u}=\vect(\mathcal{X})$.
\end{enumerate}
The computational cost of the algorithm based on Table \ref{table:2} is $\mathcal {O}(N^4)=\mathcal {O}(n^{4/3})$.  It should be mentioned here that the new formula is superior to the cyclic reduction \cite{Walter1997} and successive overrelaxation methods \cite{Young1954} that have a complexity of $\mathcal {O}(n^{3/2})$. Moreover, the new formula can be extended to higher dimensions.
\begin{table}[H]
\centering
\begin{tabular}{|c| c| c| c|} 
 \hline
 Step &  $+, -$ & $\times,  /$ & Total \\
\hline\hline
2 & $2N^3$ & $N^3$ & $3N^3$ \\ 
3 & $3N^3(N-1)$  & $3N^4$ & $6N^4-3N^3$  \\ 
4 & 0 & $N^3$ & $N^3$  \\ 
5 & $3N^3(N-1)$  & $3N^4$ & $6N^4-3N^3$  \\ 
  \hline
 & & & $12N^4-2N^3$  \\ 
 \hline
\end{tabular}
\caption{The computational cost of the algorithm for 3D Poisson equation}
\label{table:2}
\end{table}

To study their performance, scalability and accuracy of the two algorithms for solving the $2D$ and $3D$ Poisson equations, we implemented them using MATLAB scripts. All numerical experiments were carried out using Matlab R2020a running on an AMD Ryzen 5 3600 processor with 8.0GB RAM. Figure 1 shows the execution time at various grid resolutions. Importantly, as expected the $3D$ solver demonstrates significantly better performance and scalability compared to the $2D$ solver. This suggests that the algorithm used for the $3D$ case is more efficient in handling larger grid sizes, making it desirable for high resolution simulations.

We next considered the $2D$ Poisson equation
\begin{equation}
u_{xx}+u_{yy}= -200 \pi^2 \sin(10\pi x) \sin(10\pi y), u|_{\partial \Omega}=0,
\end{equation}
on the domain $[-1,1]\times[-1,1]$, which has the exact solution 
\begin{equation}
u(x,y)=\sin(10\pi x) \sin(10\pi y),
\end{equation}
to demonstrate the accuracy of the $2D$ solver. Figure 2 presents the numerical solution at various grid resolutions (column 1) alongside the exact solution (column 2) and the error difference between them (column 3). It is evident that the numerical solution progressively improves as the resolution of the grid increases.
\begin{figure}[H]
\centering
\includegraphics[totalheight=8cm]{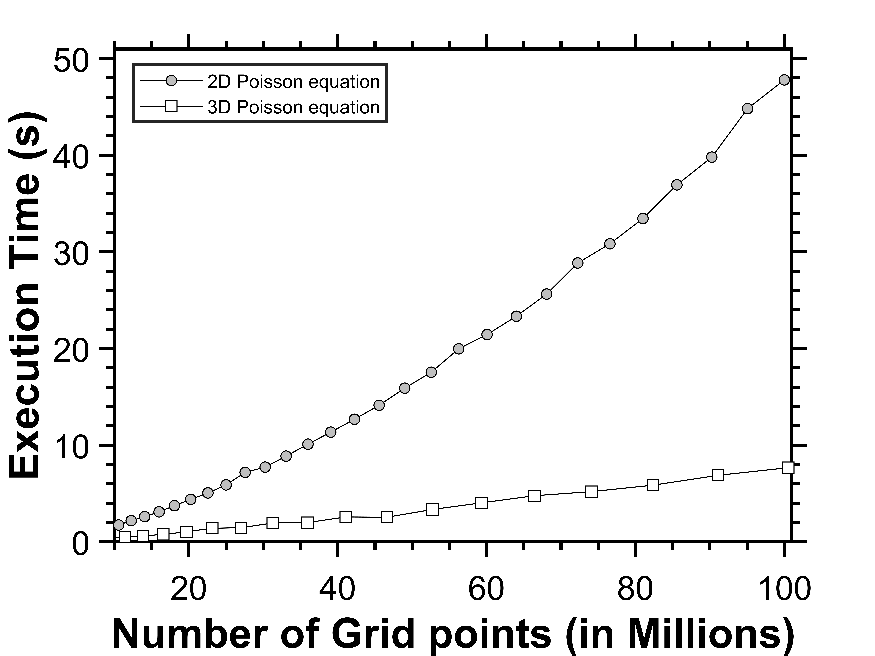}
\caption{Performance analysis showcasing the execution time for solving $2D$ and $3D$ Poisson equations at various grid resolutions.}
\label{Figure 1}
\end{figure}

\subsection{Solution of the 2D steady state convection-diffusion equation}
The 2D steady state convection diffusion equation defined on a domain $\Omega=[0,1]\times [0,1] $ 
\begin{equation}
-\nu \nabla^2 u+c^T\nabla u=f, u|_{\partial \Omega}=0,
\end{equation}
where $\nu$ is the diffusion coefficient and $c_1, c_2$ are the dominant convection coefficients, each chosen to be 1.

The right-hand side $f$ is selected such that the solution is
\begin{equation}
u\left(x,y\right)=16(x-x^2)(y-y^2),
\end{equation}
As mentiond in \cite{Grasedyck2004}, a standard finite difference discretization on a uniform grid is applied for the diffusion term, while a second-order convergent scheme (Fromm's scheme) is used for the convection term. The resulting discrete system forms a linear system of the form
\begin{equation}\label{2D_convection_diffusion}
(\mathbf{A}_1 \oplus \mathbf{A}_2)\mathbf{u}=\mathbf{b},\\
\end{equation}
where $\mathbf{A}_i, i=1,2,$ is given by
\begin{equation}
\mathbf{A}_i= \frac{\nu}{h^2}\begin{pmatrix}
2& -1 & 0 & \cdots & 0 \\
-1  & 2 &-1 & \ddots & \vdots \\
0 & -1 & 2 & \ddots & 0\\
\vdots & \ddots & \ddots & \ddots & -1 \\
0 & \cdots & 0 & -1  & 2
\end{pmatrix}_{N\times N}
+\frac{c_i}{4h}
\begin{pmatrix}
3& -5 & 1 & \cdots & 0 \\
1  & 3 &-5 & \ddots & \vdots \\
0 & 1 & 3 & \ddots & 1\\
\vdots & \ddots & \ddots & \ddots & -5 \\
0 & \cdots & 0 & 1  & 3
\end{pmatrix}_{N\times N}
,~i=1,2.
\end{equation}
If we define the $N \times N$ matrices $X$ and $Y$ such that $\mathbf{u}=\vect(\mathbf{X})$ and $\mathbf{b}=\vect(\mathbf{Y})$, then the solution of the system (\ref{2D_convection_diffusion}) can be given using Corollary (\ref{co2})  as follows:
\begin{enumerate}
\item[Step 1:]  Compute the upper triangular matrices $\mathbf{T}^{(i)}$ and unitary matrices $\mathbf{U}^{(i)}$ using the Schur decomposition $U^{(i)} \mathbf{T}^{(i)} \mathbf{U}^{(i)^*}$ of $\mathbf{A}^{(i)}$ for $i=1,2$.
\item[Step 2:] Compute the matrices  $\tilde{\mathbf{T}^{(i)}}=\mathbf{T}^{(i)}-diag(\mathbf{T}^{(i)})$ for $i=1, 2$.
\item[Step 3:]  Compute the matrix $\mathbf{C}$ such that $C_{ij}=\frac{1}{\mathbf{\Lambda}^{(1)}_i+\mathbf{\Lambda}^{(2)}_j}$.
\item[Step 4:] Compute the matrix $\mathbf{K}_0=\mathbf{C} \odot  \left( \mathbf{U}^{(1)^*}\mathbf{Y}  \mathbf{U}^{(2)}\right)$.
\item[Step 5:] Compute the matrices  $\mathbf{K}_{j+1}=\mathbf{C} \odot \left( \tilde{\mathbf{T}}^{(1)} \mathbf{K}_j +K_j \tilde{ \mathbf{T}}^{(2)^T}\right)$ for $j=0,\dots, N^2-1$.
\item[Step 6:] Compute the matrix $\mathbf{X}=\mathbf{U}^{(1)} \left(\sum\limits_{j=0}^{N^2-1}{(-1)}^j \mathbf{K}_j\right) \mathbf{U}^{(2)^*}$.
\item[Step 7:] Compute the solution $\mathbf{u}=\vect(\mathbf{X})$.
\end{enumerate}
Figure 3 illustrates the relative error $\epsilon$ between the numerical and analytical solutions of the 2D convection-diffusion equation as a function of the number of iterations for various grid resolutions. $\epsilon$ is given by
\begin{equation}
\epsilon=\frac{\lVert \mathbf{u}_{analtical}-\mathbf{u}_{numerical}\rVert}{\lVert \mathbf{u}_{analytical}\rVert}.
\end{equation}
The error decreases steadily and plateau to the discretization error level after approximately 5 iterations for all grid resolutions. This demonstrates the rapid convergence of our method. As expected, finer grids (e.g.,$4000 \times 4000$) exhibit significantly lower discretization error compared to coarser grids. Figure 3 also suggests that the convergence rate of the linear system solution is independent of size as indicated by the similar slope in the decreasing regions of the error curves across different grid resolutions.

Importantly, Remark 4.6 and Figure 3 suggest the solution can be obtained after small number of iterations. Consequently, the computational complexity can be approximated as $\mathcal {O}(n^{3/2})$.
\begin{figure}[H]
\centering
\includegraphics[totalheight=8cm]{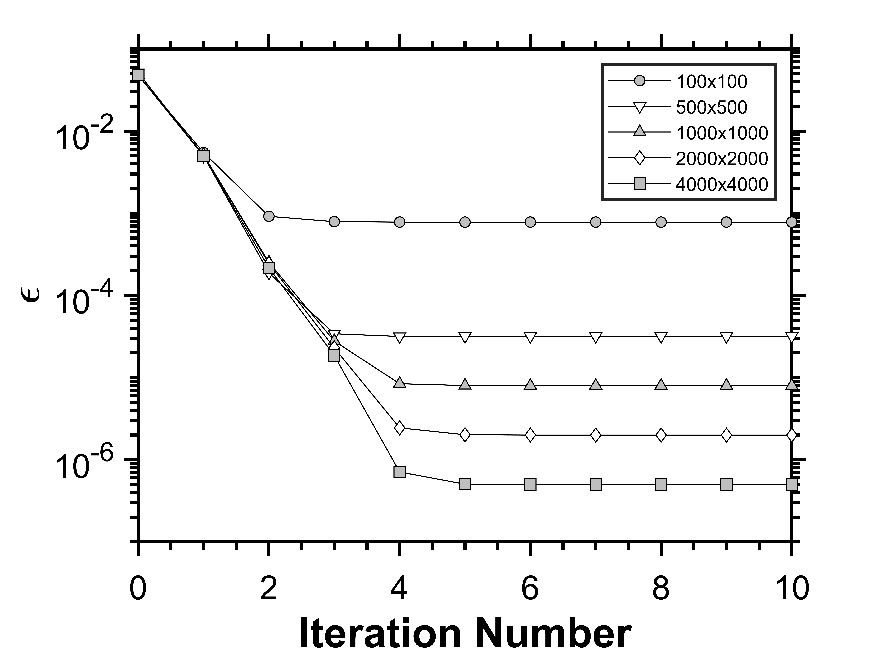}
\caption{The relative error between the numerical and analytical solutions of the 2D convection-diffusion equation as a function of the number of iterations for various grid resolutions.}
\label{Figure 3}
\end{figure}
\section{Conclusion}
The linear system involving the Kronecker sum of matrices is related to the matrix and tensor forms of Sylvester equation. Consequently, the new tensor-based formula can be seen as a novel direct solver for both the matrix and tensor forms of the Sylvester equation. Since the Sylvester equation has a unique solution for $\mathbf{X}$ when  $\lambda_1+\lambda_2 \ne 0,~\forall \lambda_1 \in \sigma(\mathbf{A}),\lambda_2 \in \sigma(\mathbf{B})$, the tensor-based formula shows that a Sylvester tensor equation has a unique solution for $\mathcal{X}$ when  $\lambda_1+\lambda_2+\lambda_3 \ne 0,~\forall \lambda_1 \in \sigma(\mathbf{A}),\lambda_2 \in \sigma(\mathbf{B}), \lambda_3 \in \sigma(\mathbf{C})$.

The tensor-based formula is efficient when applied to solve the discrete Poisson and steady state convection diffusion equations with Dirichlet boundary conditions, outperforming cyclic reduction and successive overrelaxation methods.


Conflict of Interest: The authors declare that they have no conflict of interest.

\end{document}